# BIASED RANDOM-TO-TOP SHUFFLING

## By Johan Jonasson[1]

### *Chalmers University of Technology*

Recently Wilson [*Ann. Appl. Probab.* **14** (2004) 274–325] introduced an important new technique for lower bounding the mixing time of a Markov chain. In this paper we extend Wilson's technique to find lower bounds of the correct order for card shuffling Markov chains where at each time step a random card is picked and put at the top of the deck. Two classes of such shuffles are addressed, one where the probability that a given card is picked at a given time step depends on its *identity*, the so-called move-to-front scheme, and one where it depends on its *position*.

For the move-to-front scheme, a test function that is a combination of several different eigenvectors of the transition matrix is used. A general method for finding and using such a test function, under a natural negative dependence condition, is introduced. It is shown that the correct order of the mixing time is given by the biased coupon collector's problem corresponding to the move-to-front scheme at hand.

For the second class, a version of Wilson's technique for complex-valued eigenvalues/eigenvectors is used. Such variants were presented in [*Random Walks and Geometry* (2004) 515–532] and [*Electron. Comm. Probab.* **8** (2003) 77–85]. Here we present another such variant which seems to be the most natural one for this particular class of problems. To find the eigenvalues for the general case of the second class of problems is difficult, so we restrict attention to two special cases. In the first case the card that is moved to the top is picked uniformly at random from the bottom $k = k(n) = o(n)$ cards, and we find the lower bound $(n^3/(4\pi^2k(k-1)))\log n$. Via a coupling, an upper bound exceeding this by only a factor 4 is found. This generalizes Wilson's [*Electron. Comm. Probab.* **8** (2003) 77–85] result on the Rudvalis shuffle and Goel's [*Ann. Appl. Probab.* **16** (2006) 30–55] result on top-to-bottom shuffles. In the second case the card moved to the top is, with probability $1/2$, the bottom card and with probability $1/2$, the card at position $n - k$. Here the lower bound is again

Received April 2005; revised November 2005.

[1]Supported in part by the Swedish Research Council.

of order $(n^3/k^2)\log n$, but in this case this does not seem to be tight unless $k = O(1)$. What the correct order of mixing is in this case is an open question. We show that when $k = n/2$, it is at least $\Theta(n^2)$.

**1. Introduction.** How many steps does a Markov chain need to get close to stationarity? This question has attracted a great amount of interest during the last few decades, partly because computer development has supplied the possibility of powerful MCMC simulation techniques.

Much interest has been focused on *card shuffling* chains, that is, chains whose state space is the symmetric group, $S_n$. During the 1980s and early 1990s a lot of progress was made and good upper and lower bounds, and often cutoffs, were found for a great variety of different card shuffling techniques. Among these, the Bayer–Diaconis $\frac{3}{2}\log_2 n$ cutoff for the riffle shuffle (see [2]) is the most celebrated result. In the mid-90s the subject fell into a relative silence as the available techniques did not suffice to solve the remaining unsolved problems.

However, quite recently the subject was revitalized by Wilson's [12] introduction of a powerful new technique to lower bound the mixing time. The idea is to find an easily expressed (right) eigenvector of the transition matrix, corresponding to an eigenvalue close to 1, and use this eigenvector to construct a test function. In [12] Wilson established essentially tight (i.e., tight up to a constant) lower bounds for neighbor-transposing shuffling and also for so-called lozenge tilings. A variant of Wilson's technique for complex-valued eigenvalues/eigenvectors was used by Mossel, Peres and Sinclair [8] to establish a tight lower bound for the cyclic-to-random shuffle. Wilson himself used another such variant to establish the $\Theta(n^3\log n)$ mixing time for the Rudvalis shuffle. (This result was generalized by Goel [4] who considered random-to-bottom shuffling in general.) Jonasson [6] developed a version of Wilson's technique where the test function is only needed to be constructed from something sufficiently close to an eigenvector, thereby being able to establish the $\Theta(n^2\log n)$ mixing time for the overhand shuffle.

In the present paper this development is continued. We consider random-to-top shuffles, that is, card shuffling chains where every transition is such that some card is moved from its present position to the top of the deck without affecting the relative positions of the other cards. The shuffles will be biased, that is, such that the card taken to the top is in general not chosen uniformly at random from all the cards. We consider two classes of biased random-to-top shuffles:

(1) Each *card*, $k$, is once and for all given a probability $p_k$ and the card taken to the top is chosen according to these probabilities, regardless of the present order of the deck.

(2) Same as (1), except that the probability $p_k$ is assigned to *position $k$* in the deck.

The case described by (1) is often referred to as the *move-to-front scheme*. It will always be assumed (without loss of generality) that all the $p_k$'s are positive. The move-to-front scheme clearly describes an irreducible aperiodic Markov chain, so it converges to its stationary distribution whose nature is such that the probability of having card $c_j$ in position $j$, $j \in [n]$, equals

$$p_{c_1} \cdot \frac{p_{c_2}}{1 - p_{c_1}} \cdot \frac{p_{c_3}}{1 - p_{c_1} - p_{c_2}} \cdot \cdots \cdot \frac{p_{c_n}}{1 - \sum_{j=1}^{n-1} p_{c_j}};$$

see, for example, [9]. To the best of my knowledge, convergence rates for the move-to-front scheme have only been studied for the case where all the $p_k$'s are equal; the ordinary random-to-top shuffle for which one has a cutoff at $n \log n$ shuffles; see [1]. Rodrigues [9] studied a variant of the move-to-front scheme where the identities of the cards moved to top at different shuffles are dependent. For the ordinary top-to-random shuffle, the mixing time is given by the classical coupon-collector's problem. In Section 3 we show that the corresponding biased coupon-collector's problem in the general case gives the correct answer up to a constant factor. Our lower bounds are established via an extension of Wilson's technique where several different eigenvalues/eigenvectors are used to construct the test function.

Case (2) with $p_{n-1} = p_n = 1/2$ is known as the Rudvalis shuffle. Therefore, we refer to case (2) from now on as *generalized Rudvalis* (or GR) shuffles. The special case when, for some $k \leq n$, the card moved to the top is chosen uniformly among the bottom $k = k(n)$ cards is, in the terminology of Goel [4], referred to as the *bottom-to-top* shuffle. When $k = O(1)$, it is known that the mixing time is of order $n^3 \log n$; for the Rudvalis shuffle, see [5] and [11] and for the general case, see [4]. Goel also estimated the mixing time for other values of $k$: When $k = \Theta(n)$, it is shown that $k$ is at least of order $n \log n$ and at most of order $n^2 \log n$ and that when $k$ is sufficiently close to $n$, the correct order is $n \log n$. In Section 4, which is devoted to GR shuffles, we improve the previous results on the bottom-to-top shuffles. For $k = o(n)$, we find that the correct order of the mixing time is $(n^3/k^2) \log n$ and find upper and lower bounds that only differ by a factor 4. In the case $k = \Theta(n)$ it is shown that the mixing time is $\Theta(n \log n)$. The upper bounds are found via coupling and to estimate the coupling time, we use the trick of imagining the top of the deck as shifting cyclically one step down the deck for every time step. Then if a card is only observed at times it is touched, its motion is described by a random walk whose step size has mean 0 and variance of order $k^2$. This can then be used together with the nature of the coupling to conclude the time to get a particular card matched is roughly bounded by $n$ times the time taken for a random walk of this type started at $n/2$ to exit the interval $[0, n]$. The lower bounds are found using a version of Wilson's technique designed to handle complex-valued eigenvalues/eigenvectors. Wilson [11] and Saloff-Coste [10] designed other such variants. They also pointed out that if one tries to apply Wilson's original lemma in complex-valued cases without further thought, one often ends up with results that are much too weak to be of interest. Our version is the perhaps most straightforward adaption of Wilson's technique to the complex-valued case and it does not need the chain under study to be lifted to a larger state space. It works well for the present applications and we believe that it is fairly generally applicable, but it would, for example, not solve all problems in [11].

The lower bound technique, in principle, works for GR shuffles in general. However, in practice, it is very difficult to find the eigenvalues, even for the motion of a single card, but in the simplest cases. In addition to the bottom-to-top shuffle, we also carry out the calculations for the cases $p_{n-k} = p_n = 1/2$ for different $k$'s and find a lower bound which is almost identical to the one for the bottom-to-top shuffle with the same $k$. However, unlike the bottom-to-top case, this lower bound does not seem to be of the correct order in this case, unless $k = O(1)$. We show that when $k = n/2$, the mixing time is at least $\Theta(n^2)$. In fact, single-card motion takes $\Theta(n^2)$ steps to mix, so this is an example of a situation where the second eigenvalue for single-card motion does not capture the correct order of the mixing time, not even for the single-card chain itself. Unfortunately we have not been able to produce any good upper bound, and it is a wide open question what the order of mixing really is.

The next section gives the necessary preliminaries.

## 2. Preliminaries.

2.1. *Basic definitions.* The most common way to measure the distance between two probability measures $\mu$ and $\nu$ on a finite set $S$ is by the *total variation norm* given by

$$\|\mu - \nu\| := \tfrac{1}{2} \sum_{s \in S} |\mu(s) - \nu(s)| = \max_{A \subseteq S} (\mu(A) - \nu(A)).$$

If $\{X_t\}_{t=0}^{\infty}$ is an aperiodic irreducible Markov chain on state space $S$ and with stationary distribution $\pi$, started from a fixed state $s$, then its *mixing time*, $\tau_{\mathrm{mix}}$, is defined via

$$\tau(s) := \min\{t : \|P(X_t \in \cdot) - \pi\| \le \tfrac{1}{4}\}$$

and

$$\tau_{\mathrm{mix}} := \max_s \tau(s).$$

The mixing time is often expressed in terms of some measure of the size of the state space. When doing so one implicitly considers a sequence of Markov chains $\{X_t^n\}$, $n = 1, 2, 3, \ldots$, on state spaces $S^n$ and with stationary distributions $\pi^n$, where the state spaces are such that $|S^n| \uparrow \infty$ in some natural way. In our case $S^n = S_n$, the symmetric group on $n$ cards, and the mixing time is expressed in relation to the number of cards, $n$, as $n \to \infty$. The sequence of Markov chains is said to have a *cutoff* at $\tau_{\text{mix}} := \tau_{\text{mix}}^n$ if, for every $a > 0$,

$$\lim_{n \to \infty} \|P(X_{(1+a)\tau_{\text{mix}}}^n \in \cdot) - \pi^n\| = 0$$

and

$$\lim_{n \to \infty} \|P(X_{(1-a)\tau_{\text{mix}}}^n \in \cdot) - \pi^n\| = 1.$$

2.2. *Coupling.* A common technique to find upper bounds on the mixing time is *coupling*: Suppose $\{X_t\}$ is the Markov chain under study, started from a fixed state $s$ maximizing $\tau(s)$, and that $\{Y_t\}$ is a chain with the same transition rule, but started from stationarity. Suppose also that the updates of the two chains have been set up, or *coupled*, in such a way that as soon as $X_{t_0} = Y_{t_0}$ for some $t_0$, then $X_t = Y_t$ for all $t \geq t_0$. Put

$$T := \min\{t : X_t = Y_t\}.$$

Then $T$ is called the *coupling time* and the *coupling inequality* (see, e.g., [7], Section I.2) states that, for all $t$,

$$\|P(X_t \in \cdot) - \pi\| \leq P(T > t).$$

The coupling inequality follows from the definition of total variation norm and the simple observation that, on the event $\{T \leq t\}$, $X_t$ and $Y_t$ are the same.

2.3. *Wilson's technique.* Here we state and prove a basic version of Wilson's technique that will later develop in a few different ways. Let $\{X_t\}_{t=0}^\infty$ be an irreducible aperiodic Markov chain on state space $S$ and with stationary distribution $\pi$, starting from a fixed state $s_0$. Let $\Phi : S \to \mathbb{C}$ and $\gamma \in (0, 1/2)$ be such that, for every $s \in S$,

$$\mathbb{E}[\Phi(X_{t+1})|X_t = s] = (1 - \gamma)\Phi(X_t),$$

that is, $1 - \gamma$ is an eigenvalue for the transition matrix of the chain and $\Phi$ is a corresponding eigenvector. Let $R$ be a number such that

$$R \geq \max_{s \in S} \mathbb{E}[|\Phi(X_{t+1}) - \Phi(X_t)|^2 | X_t = s].$$

THEOREM 2.1. *Fix $a > 0$ and let*

$$T := \frac{\log |\Phi(s_0)| - (1/2)\log 4R/(\gamma a)}{-\log(1 - \gamma)}.$$

*Then $\|P(X_t \in \cdot) - \pi\| \geq 1 - a$ for all $t \leq T$.*

PROOF. By induction, $\mathbb{E}\Phi(X_t) = (1-\gamma)^t \Phi(s_0)$. Therefore, $\mathbb{E}\Phi(X_\infty) = 0$, with $X_\infty$ denoting a state chosen from the stationary distribution. Put $\Delta\Phi = \Phi(X_{t+1}) - \Phi(X_t)$. We have that

$$\mathbb{E}[|\Phi(X_{t+1})|^2 | X_t] = \mathbb{E}[|\Phi(X_{t+1}) + \Delta\Phi|^2 | X_t]$$
$$= (1-2\gamma)|\Phi(X_t)|^2 + \mathbb{E}[|\Delta\Phi|^2 | X_t]$$
$$\leq (1-2\gamma)|\Phi(X_t)|^2 + R.$$

By induction using that $\gamma \leq 1/2$,

$$\mathbb{E}|\Phi(X_t)|^2 \leq (1-2\gamma)^t |\Phi(s_0)|^2 + \frac{R}{2\gamma}.$$

Thus, again using that $\gamma \leq 1/2$,

$$\mathrm{Var}\,\Phi(X_t) = \mathbb{E}|\Phi(X_t)|^2 - |\mathbb{E}\Phi(X_t)|^2$$
$$= ((1-2\gamma)^t - (1-\gamma)^{2t})|\Phi(s_0)|^2 + \frac{R}{2\gamma} \leq \frac{R}{2\gamma}.$$

By Chebyshev's inequality,

$$P\left(|\Phi(X_t) - \mathbb{E}\Phi(X_t)| \geq \sqrt{\frac{R}{\gamma a}}\right) \leq \frac{a}{2}$$

and, on letting $t \to \infty$,

$$P\left(|\Phi(X_\infty)| \geq \sqrt{\frac{R}{\gamma a}}\right) \leq \frac{a}{2}.$$

Thus, if $t$ is such that $|\mathbb{E}\Phi(X_t)| \geq 2\sqrt{R/(\gamma a)}$, we have $\|P(X_t \in \cdot) - \pi\| \geq 1 - a$. This, however, is the case precisely when $t \leq T$. $\quad\square$

**3. The move-to-front scheme.** Recall the move-to-front scheme. A deck of $n$ cards is shuffled by the following rule: To each card $k$, $k \in [n]$, we attach once and for all a probability $p_k$. We assume without loss of generality that $p_k > 0$ for every $k$ and that $p_1 \geq p_2 \geq \cdots \geq p_n$. We will also for technical reasons work under the mild condition that $p_k \leq 1/3$ for every $k$. At each time step a card is chosen according to these probabilities and then moved to the top of the deck without altering the relative positions among the other cards. Let $\{X_t\}_{t=0}^\infty$ denote the Markov chain on the symmetric group $S_n$ defined in this way, started from a permutation $X_0 = s_0$ maximizing the time taken to reach stationarity.

3.1. *Upper bound.* To find an upper bound on the mixing time, we use coupling. For this we need another deck $\{Y_t\}$ started from stationarity to be updated according to the same transition rule. The coupling is given by the following simple rule: At each time step let the same card be moved to the top for the two decks. Let $T = \min\{t : X_t = Y_t\}$ be the coupling time. Considering only one of the decks, a moment's thought reveals that at the first time that all but one card has at least once been picked to the top, the order of the cards no longer depends on their starting positions, but can be read off completely from in what order they have been moved to the top. More precisely, the later a card has been picked, the higher up in the deck it is, with the still untouched card at the very bottom of the deck. Now considering again both decks, this means that at this time the two decks must agree. Therefore, the problem can be solved by solving the biased coupon collector's problem naturally appearing from the different card probabilities:

$$P(T \geq t) \leq P(\text{at least two cards not touched by time } t) \leq \sum_{k=1}^{n-1} (1 - p_k)^t.$$

Thus, putting

$$\tau_u := \min\left\{ t : \sum_{k=1}^{n-1} (1 - p_k)^t \leq \tfrac{1}{4} \right\},$$

the coupling inequality tells us that $\tau_{\mathrm{mix}} \leq \tau_u$. Our next task is to prove that $\tau_u$ is, up to a constant, the correct mixing time.

3.2. *Lower bound.* Because of the fact that different cards are picked with different probabilities, it is difficult to find eigenvectors for single eigenvalues that, evaluated at $X_t$, are sufficiently concentrated around their mean to provide a good enough test function for sharp lower bounds. We shall circumvent this problem by combining eigenvectors corresponding to several different eigenvalues. For a general treatment of the idea, we assume that the setting is the same as in Section 2.3, but that $(1 - \gamma_j, \Phi_j)$, $j = 1, \ldots, m$, are possibly different eigenvalue/eigenvector pairs for the transition matrix. For simplicity, assume that the $\Phi_j$'s are scaled in such a way that, for all $j$, $\max_{s \in S} \mathbb{E}[|\Delta \Phi_j|^2 | X_t = s] \leq \gamma_j$. [The scaling does not in any way affect the lower bound that comes out in the end, what matters is the relation between the variance of $\Phi_j(X_t)$ and $\Phi_j(s_0)$.] Assume also that, for every $j$, $k$ and $t$,

$$\mathrm{Cov}(\Phi_j(X_t), \Phi_k(X_t)) \leq 0.$$

The test function we shall use is the following:

$$\Phi(X_t) = \Phi^t(X_t) := \sum_{j=1}^{m} a_j \Phi_j(X_t),$$

where the $a_j = a_j(t)$, $j \in [m]$, form a unit vector, that is, $\sum_{j=1}^m a_j^2 = 1$. The precise choice of the $a_j$'s should be made in such a way that the lower bound achieved is optimized. Reasoning as in the proof of Theorem 2.1, $\operatorname{Var} \Phi_j(X_t) \leq 1/2$. Thus, by the negative covariance assumption, $\operatorname{Var} \Phi(X_t) \leq 1/2$ and by letting $t \to \infty$, $\operatorname{Var} \Phi(X_\infty) \leq 1/2$. Note also that $\mathbb{E}\Phi(X_t) = \sum_{j=1}^m a_j(1 - \gamma_j)^t \Phi_j(s_0)$ and $\mathbb{E}\Phi(X_\infty) = 0$. Continuing along the lines of the proof of Theorem 2.1 now reveals that $\|P(X_t \in \cdot) - \pi\| \geq 1/3$ as long as $t$ is such that $|\mathbb{E}\Phi(X_t)| \geq \sqrt{6}$, that is, if $|\sum_{j=1}^m a_j(1-\gamma_j)^t \Phi_j(s_0)| \geq \sqrt{6}$ for an optimal choice of the $a_j$'s.

Let us now move back to the move-to-front scheme. Let the deck start with card $n$ in position 1, card $n-1$ in position 2,..., card 1 in position $n$ and denote this state as before by $s_0$. Assume for simplicity that $n$ is even. For every odd $j \in [n]$, define

$$\Phi_j(X_t) = \begin{cases} \dfrac{p_j}{p_j + p_{j+1}}, & \text{if } X_t(j+1) \leq X_t(j), \\ -\dfrac{p_{j+1}}{p_j + p_{j+1}}, & \text{if } X_t(j) < X_t(j+1). \end{cases}$$

In words, $\Phi_j$ is assigned the positive value when card $j+1$ is above card $j$ in the deck and the negative value when the opposite situation occurs. Then $\Phi_j$ is an eigenvector for the eigenvalue $1 - \gamma_j$, where $\gamma_j := p_j + p_{j+1}$. We have that $\mathbb{E}[|\Delta\Phi_j|^2|X_t] \leq p_j \leq \gamma_j$ so the general situation just described applies if it can be shown that the $\Phi_j$'s are negatively correlated. However, the simplicity of the situation allows for a direct calculation of the covariance of $\Phi_j(X_t)$ and $\Phi_k(X_t)$, $j \neq k$:

$$\mathbb{E}[\Phi_j(X_t)\Phi_k(X_t)] = (1 - \gamma_j - \gamma_k)^t \frac{p_j}{\gamma_j} \frac{p_k}{\gamma_k}$$

and

$$\mathbb{E}[\Phi_j(X_t)] = (1 - \gamma_j)^t \frac{p_j}{\gamma_j},$$

so

$$\operatorname{Cov}(\Phi_j(X_t), \Phi_k(X_t)) = ((1 - \gamma_j - \gamma_k)^t - (1 - \gamma_j)^t(1 - \gamma_k)^t)\frac{p_j p_k}{\gamma_j \gamma_k} \leq 0.$$

Since $\Phi_j(s_0) = p_j/\gamma_j \geq 1/2$, it follows from above that the variation distance from stationarity is at least $1/3$ when $t$ is such that $\sum_{j:j \text{ odd}} a_j(1 - \gamma_j)^t \geq 2\sqrt{6}$. In order to make this bound on $t$ as good as possible, we pick the $a_j$'s so that the left-hand side is maximized, that is, $a_j = (1-\gamma_j)^t/(\sum_{k:k \text{ odd}}(1-\gamma_k)^{2t})^{1/2}$. Doing so we get that variation distance is at least $1/3$ for $t$ such that

$$\sum_{j:j \text{ odd}} (1 - \gamma_j)^{2t} \geq 24.$$

In order to put this in relation to $\tau_u$ above, recall that all $p_j$'s are assumed not to exceed $1/3$. Therefore, $(1 - \gamma_j)^{2t} \geq (1 - 2p_j)^{2t} \geq (1 - p_j)^{6t}$. Note also that $\sum_{j\,:\,j\text{ odd}}(1 - p_j)^{6t} \geq \frac{1}{2}\sum_{j=1}^{n-1}(1 - p_j)^{6t}$. Therefore, variation distance is at least $1/3$ as long as

$$\sum_{j=1}^{n-1}(1 - p_j)^{6t} \geq 48.$$

Put $\tau_0$ for the largest $t$ for which this inequality holds. Now in case $\tau_0 < \log(3/4)/\log(1 - p_{n-1})$, then $(1 - p_{n-1})^{\tau_0} > 3/4$. In this case $\tau_0$ is not a good lower bound and a better one is given by the trivial one $\tau_1 := \max\{t\,:\,(1 - p_{n-1})^t > 3/4\}$. (Trivial because the probability at stationarity of having card $n-1$ above card $p_n$ is at least $1/2$ and it is less than $1/4$ at time $\tau_1$.) However, since $\sum_{j=1}^{n-1}(1 - p_j)^{6(\tau_1+1)} < 48$ and $(1 - p_j)^{\tau_1+1} \leq 3/4$ for all $j \leq n-1$,

$$\sum_{j=1}^{n-1}(1 - p_j)^{25(\tau_1+1)} \leq 48 \cdot \left(\tfrac{3}{4}\right)^{19} < \tfrac{1}{4}.$$

Hence, $\tau_u \leq 25(\tau_1 + 1)$ and so the biased coupon collector's problem yields the correct mixing time up to a factor of 25. Assume now that $(1 - p_{n-1})^{\tau_0} \leq 3/4$ so that $(1 - p_j)^{\tau_0+1} < 3/4$ for all $j$. Then

$$(1 - p_j)^{25(\tau_0+1)} < \left(\tfrac{3}{4}\right)^{19}(1 - p_j)^{6(\tau_0+1)} < \tfrac{1}{184}(1 - p_j)^{6(\tau_0+1)}$$

and so

$$\sum_{j=1}^{n-1}(1 - p_j)^{25(\tau_0+1)} < \frac{1}{184}\sum_{j=1}^{n-1}(1 - p_j)^{6(\tau_0+1)} \leq \tfrac{1}{4}.$$

Thus, $\tau_u \leq 25(\tau_0 + 1)$. The following theorem summarizes the results of this section.

THEOREM 3.1. *Consider the move-to-front scheme with card probabilities $p_j$, $j \in [n]$, with $p_j \leq 1/3$ for every $j$. Put*

$$\tau_u := \min\left\{t\,:\,\sum_{j=1}^{n-1}(1 - p_j)^t \leq \tfrac{1}{4}\right\}.$$

*Then*

$$\tfrac{1}{25}\tau_u - 1 \leq \tau_{\text{mix}} \leq \tau_u.$$

EXAMPLES. (a) The ordinary random-to-top shuffle. Here $p_j = 1/n$ for every $j$. We get $\tau_u = (1 + o(1))n\log n$ and $\frac{1}{25}n\log n \leq \tau_{\text{mix}} \leq n\log n$. Of course, it is well known that $n\log n$ is a cutoff in this case.

(b) Put

$$p_j = \begin{cases} \dfrac{2}{n+1}, & 1 \leq j \leq \dfrac{n}{2}, \\ \dfrac{2}{n(n+1)}, & \dfrac{n}{2} + 1 \leq j \leq n. \end{cases}$$

Here $\tau_u = (1 + o(1))\frac{1}{2}n^2 \log n$, so the correct order of the mixing time is $n^2 \log n$.

(c) Let $p_j = j^{-1}/(\sum_{k=1}^n k^{-1})$. Then with $t = cn(\log n)^2$, it is readily seen that

$$\sum_{j=1}^{n-1} (1 - p_j)^t = \begin{cases} \Omega(1), & \text{when } c < 1, \\ o(1), & \text{when } c > 1. \end{cases}$$

Therefore, $\tau_u = (1 + o(1))n(\log n)^2$ and the order of the mixing time is $n(\log n)^2$.

(d) Put $p_j = 2(n+1-j)/(n(n+1))$. With $t = cn^2$,

$$\sum_{j=1}^{n-1} (1 - p_j)^t \leq \sum_{j=1}^{n-1} e^{-2cj} \leq \frac{e^{-2c}}{1 - e^{-2c}} \leq \frac{1}{4}$$

if, say, $c \geq 1$. Therefore, $n^2$ is the correct order of the mixing time. This is thus a case where the time taken to reach stationarity depends almost entirely on the time taken to find the few most unprobable cards.

## 4. Generalized Rudvalis shuffles.

4.1. *Lower bounds.* The eigenvalues for the GR shuffles that are at least reasonably accessible are those for the motion of a single card. These typically turn out to be complex. Now taking a closer look at the proof of Theorem 2.1 reveals that there is nothing in it that technically prevents $\lambda := 1 - \gamma$ from being complex-valued. However, the typical situation is such that $\gamma$ is of much larger order than $1 - |\lambda|$, but it is the latter that indicates the correct mixing time. Therefore, one would like a variant of Theorem 2.1, where one in effect works with $1 - |\lambda|$ rather than $\gamma$. Wilson [11] developed such a method in order to take care of the original Rudvalis shuffle and some variants of it. One ingredient in his method is to extend the Markov chain to a larger state space by incorporating time into it. In our case such an extension seems unnecessary and the following variant is the most natural way to attack the problems at hand:

Let the setting be as in Theorem 2.1, with the exception that the eigenvalue is now $(1 - \gamma)e^{i\theta}$ for some $\gamma \in (0, 1/2]$ and some $\theta \in [0, \pi]$ and $R$ is such that

$$R \geq \max_{s \in S} \mathbb{E}[|e^{-i\theta}\Phi(X_{t+1}) - \Phi(X_t)|^2 | X_t = s].$$

THEOREM 4.1. *Assume the setting above, pick $a > 0$ and let $T$ be as in Theorem 2.1. Then for $t \le T$,*

$$\|P(X_t \in \cdot) - \pi\| \ge 1 - a.$$

PROOF. As in the proof of Theorem 2.1, $\mathbb{E}\Phi(X_t) = (1 - \gamma)^t e^{it\theta} \Phi(s_0)$ and $\mathbb{E}\Phi(X_\infty) = 0$. Put, for $t = 0, 1, 2, 3, \ldots$,

$$\Psi_t(X_t) = e^{-it\theta} \Phi(X_t).$$

Then for every $t$, $\mathbb{E}\Psi_t(X_t) = (1 - \gamma)^t \Phi(s_0)$ and $\mathbb{E}\Psi_t(X_\infty) = 0$ and

$$R \ge \max_{s \in S} \mathbb{E}[|\Psi_{t+1}(X_{t+1}) - \Psi_t(X_t)|^2 | X_t = s].$$

Now completely mimic the proof of Theorem 2.1 with $\Psi_t(X_t)$ playing the rôle of $\Phi(X_t)$. Doing so yields the desired result. □

Recall that in the general GR shuffle setting, the card moved to the top is chosen from position $k$ with probability $p_k$, $k \in [n]$. This entails that the motion of a single card, $c$, is in itself a Markov chain: Given that $c$ is in position $k$ at time $t$, then at time $t + 1$, $c$ is in position 1 with probability $p_1$, in position $k$ with probability $m_k := \sum_{j<k} p_j$ and in position $k + 1$ with probability $M_k := \sum_{j>k} p_j$. Putting $\lambda$ for an eigenvalue of this chain and $\mathbf{x}$ for a corresponding eigenvector, these must satisfy the equations

$$\lambda x_k = p_k x_1 + m_k x_k + M_k x_{k+1},$$

$k = 1, 2, \ldots, n$. Unfortunately this system of equations is very difficult to solve in general, so from now on we focus on some special cases. However, even so the exact eigenvalues cannot be expressed on closed form and we will consequently only be able to produce expressions that are very close to the eigenvalue. When arguing that a given expression is indeed close to an eigenvalue, the following lemma is very useful. The lemma should be more or less well known, but we have not been able to find any specific reference, so a proof is supplied.

LEMMA 4.2. *Let $\mathbb{D}$ be the closed unit disc in the complex plane. Assume that $f : \mathbb{D} \to \mathbb{C}$ is analytic, $f(0) = 0$ and $|f'(z)| \ge 1$ for every $z$ in $\mathbb{D}$. Then there exists a point $z_0 \in \mathbb{D}$ such that $f(z_0) = 1$.*

PROOF. Let $V := f(\mathbb{D})$ and let $u$ be the leftmost point of the positive real axis that intersects the boundary of $V$; formally,

$$u := \inf\{a \in [0, \infty) : a \notin V\}.$$

We want to prove that $f^{-1}[0, u]$ contains a smooth path from 0 to $\partial\mathbb{D}$. The set $f^{-1}[0, u]$ may not be connected, but put $\rho$ for the component of $f^{-1}[0, u]$ that contains the origin. We claim that $\rho$ is a path of the desired type. By the assumption on $f'$, the set $\rho$ can be covered by open balls on which the restriction of $f$ is an open bijective map with analytic inverse. Since $\rho$ is compact, this cover can be taken to be finite; put $U_1, U_2, \ldots, U_n$ for the covering balls. Since $f|_{U_j}$ has a well-defined analytic inverse, $\rho \cap U_j$ is a smooth path segment for every $j$. Putting the $U_j$'s together shows that $\rho$ is a smooth path and since $\rho$ is closed, it contains its endpoints. Finally, to prove that the endpoint, $w \neq 0$, of $\rho$ is on $\partial \mathbb{D}$, assume for contradiction that $w \in \mathbb{D}^0$. Then $w$ is the center of an open ball, $O$, contained in $\mathbb{D}^0$ on which $f$ has analytic inverse. Since $f(O) \subseteq V^0$, $f(O)$ is crossed by $[0, u]$ and, hence, $O$ is crossed by a segment of $\rho$, contradicting that $w$ is the endpoint.

Next we claim that $f|_\rho$ is a bijection onto $[0, u]$: If this was not the case, then by continuity, for some $j$, $f|_{\rho \cap U_j}$ would not be bijective. (We may regard $f|_\rho$ as a continuous real-valued function defined on a closed interval of $\mathbb{R}$.) This contradicts to the bijectivity of $f|_{U_j}$.

Thus, $\rho$ can be parameterized the natural way, that is, by letting $\rho(s) = f^{-1}(s)$, $s \in [0, u]$. For an arbitrary $z \in \rho$, write $z = \rho(s)$, $s \in [0, u]$, and get, by Cauchy's theorem,

$$s = f(z) = \int_{\rho[0,s]} f'(w) \, dw = \int_0^s f'(\rho(v)) \rho'(v) \, dv.$$

Hence, $f'(\rho(v)) \rho'(v) = 1$ for almost every $v$. Whence, by hypothesis,

$$u = \int_0^u |f'(\rho(v))| |\rho'(v)| \, dv \geq \int_0^u |\rho'(v)| \, dv.$$

However, the right-hand side is the length of $\rho$ and since $\rho$ goes from the origin to the boundary of $\mathbb{D}$, its length is at least 1. Thus, $u \geq 1$, that is, $f(\mathbb{D})$ contains 1, as desired.  $\square$

A slight strengthening of the proof of Lemma 4.2 proves that, in fact, $f(\mathbb{D}) \supseteq \mathbb{D}$: Replace the positive real axis $[0, \infty)$ with any ray $e^{i\theta}[0, \infty)$, $\theta \in [0, 2\pi)$ and redefine $u$ as $u := \inf\{a \in [0, \infty) : ae^{i\theta} \notin V\}$. Then, repeating the proof with only a small adjustment for the $e^{i\theta}$-factor again gives $u \geq 1$. Since $\theta$ was arbitrary, we have shown the following:

THEOREM 4.3.   *Let* $f : \mathbb{D} \to \mathbb{C}$ *be analytic with* $f(0) = 0$ *and* $|f'(z)| \geq 1$ *for every* $z \in \mathbb{C}$. *Then* $f(\mathbb{D}) \supseteq \mathbb{D}$.

The way in which Lemma 4.2 will be used is the following: Suppose that for some $z_0$ we have that $f(z_0) = \delta$. Suppose also that $|f'(z)| \geq M$ for all $z$ within distance $\delta/M$ of $z_0$. Then by Lemma 4.2 applied to the map $z \to 1 - f(z_0 + \delta/M)$, $f$ has a zero within distance $\delta/M$ of $z_0$.

4.1.1. *The bottom-to-top shuffle.* Recall that now the card taken to the top of the deck is chosen uniformly among the bottom $k$ cards. As above, put $\lambda$ and $\mathbf{x}$ for an eigenvalue and a corresponding eigenvector for single-card motion and assume without loss of generality that $x_1 = 1$. Then the above system of equations for $(\lambda, \mathbf{x})$ becomes

$$\lambda = x_2,$$

$$\lambda x_2 = x_3,$$

$$\lambda x_3 = x_4,$$

$$\vdots$$

$$\lambda x_{n-k} = x_{n-k+1},$$

$$\lambda x_{n-k+1} = \frac{1}{k} + \frac{k-1}{k} x_{n-k+2},$$

$$\lambda x_{n-k+2} = \frac{1}{k} + \frac{1}{k} x_{n-k+2} + \frac{k-2}{k} x_{n-k+3},$$

$$\lambda x_{n-k+3} = \frac{1}{k} + \frac{2}{k} x_{n-k+3} + \frac{k-3}{k} x_{n-k+4},$$

$$\vdots$$

$$\lambda x_{n-1} = \frac{1}{k} + \frac{k-2}{k} x_{n-1} + \frac{1}{k} x_n,$$

$$\lambda x_n = \frac{1}{k} + \frac{k-1}{k} x_n.$$

Working backward, we get from the last equation that

$$x_n = \frac{1}{k\lambda - (k-1)}.$$

The second to last equation gives

$$x_{n-1} = \frac{1/k + (1/k)x_n}{\lambda - (k-2)/k} = \frac{1}{k\lambda - (k-1)} = x_n$$

or $\lambda = (k-2)/k$. Unless $\lambda = (k-2)/k$, the third to last equation now tells us that

$$x_{n-2} = \frac{1/k + (2/k)x_{n-1}}{\lambda - (k-3)/k} = \frac{1}{k\lambda - (k-1)} = x_n$$

or $\lambda = (k-3)/k$. By induction,

$$x_n = x_{n-1} = \cdots = x_{n-k+1} = \frac{1}{k\lambda - (k-1)}$$

or $\lambda \in \{0, 1/k, 2/k, \ldots, (k-2)/k\}$. By the first $n-k$ equations, $x_1 = 1$, $x_2 = \lambda$, $x_3 = \lambda^2, \ldots, x_{n-k+1} = \lambda^{n-k}$. Unless $\lambda \in \{0, 1/k, 2/k, \ldots, (k-2)/k\}$, the two expressions for $x_{n-k+1}$ must be equal. This gives the following characteristic equation for $\lambda$:

$$g(\lambda) := \lambda^{n-k+1} - \frac{k-1}{k}\lambda^{n-k} - \frac{1}{k} = 0.$$

Now assume that $k = o(n)$. Then no eigenvalue less than $1 - 2/k$ will provide a good lower bound so we can focus on $\lambda$'s solving the characteristic equation. So how do we find solutions to this equation? It is fairly easy to guess that there may be a solution close to $\lambda = e^{iw}$, where $w = 2\pi/n$, so let us first simply put $\lambda = (1-\gamma)e^{iw}$, insert this into the equation and make an approximate calculation of what $\gamma$ then should be. We have, if $\gamma$ is assumed to be very small,

$$g(\lambda) \approx (1 - n\gamma)e^{-i(k-1)w} - \frac{k-1}{k}(1 - n\gamma)e^{-ikw} - \frac{1}{k}.$$

The imaginary part of this is very small. The real part is approximately

$$(1 - n\gamma)\left(1 - \frac{(k-1)^2}{2}w^2\right) - \frac{k-1}{k}(1 - n\gamma)\left(1 - \frac{k^2w^2}{2}\right) - \frac{1}{k}$$

$$\approx -\frac{1}{k}n\gamma + \left(\frac{k(k-1)}{2} - \frac{(k-1)^2}{2}\right)w^2.$$

To make the last expression vanish, we must put

$$\gamma = \binom{k}{2}\frac{w^2}{n}.$$

We thus guess that $g(\lambda)$ has a zero very close to $\lambda_0 := (1 - \binom{k}{2}w^2/n)e^{iw}$. We now need to prove that this is indeed the case. To make things a bit easier to handle, transform the characteristic equation by putting $z = \lambda e^{-iw}$, thereby getting

$$f(z) := z^{n-k+1} - \frac{k-1}{k}e^{-iw}z^{n-k} - \frac{1}{k}e^{i(k-1)w} = 0,$$

for which we want to prove there is a root very close to $z_0 := 1 - \binom{k}{2}w^2/n$. More precisely, we will use Lemma 4.2 to show that the distance from $z_0$ to the nearest root is $O(k^3n^{-4}) = o(k^2w^2/n)$. Whence, there is a root of the form $1 - (1 + o(1))\binom{k}{2}w^2/n$. First we calculate $f'$:

$$f'(z) = z^{n-k-1}\left((n-k+1)z - \frac{(n-k)(k-1)}{k}e^{-iw}\right).$$

Next we evaluate $f'(z)$ for an arbitrary point $z$ such that $|z - z_0| < \binom{k}{2} w^2 / (2n)$. Such a point can be written as $1 - c \binom{k}{2} w^2 / n$, where $|c| \in (1/2, 3/2)$. We have

$$
\begin{aligned}
f'(z) &= (1 - o(1)) \left( (n - k + 1) \left( 1 - c \binom{k}{2} \frac{w^2}{n} \right) \right. \\
&\qquad\qquad \left. - \frac{(n-k)(k-1)}{k} (1 - O(n^{-2})) + i(-w + O(n^{-3})) \right) \\
&= (1 - o(1)) \left( \frac{k(n-k+1) - (n-k)(k-1)}{k} - O(k^2 n^{-2}) + O(n^{-1}) \right. \\
&\qquad\qquad\qquad\qquad \left. + i \left( \frac{(n-k)(k-1)}{k} w + O(n^{-2}) \right) \right) \\
&= (1 + o(1)) \left( \frac{n}{k} + i \frac{2\pi(k-1)}{k} \right).
\end{aligned}
$$

Thus, $f'(z) = (1 + o(1)) n / k$ for all $z$ within distance $\binom{k}{2} w^2 / (2n)$ of $z_0$. If it can now be shown that $f(z_0) = O(k^2 n^{-3})$, it will follow from Lemma 4.2 that $f$ has a zero within distance $O(k^3 n^{-4})$ of $z_0$, as desired. However,

$$
\begin{aligned}
f(z_0) &= z_0^{n-k+1} - \frac{k-1}{k} e^{-iw} z_0^{n-k} - \frac{1}{k} e^{i(k-1)w} \\
&= 1 - \frac{n-k+1}{n} \binom{k}{2} w^2 - \frac{k-1}{k} \left( 1 - \frac{n-k}{n} \binom{k}{2} w^2 \right) \left( 1 - \frac{w^2}{2} - iw \right) \\
&\qquad - \frac{1}{k} \left( 1 - \frac{(k-1)^2 w^2}{2} + i(k-1)w \right) + O(k^2 n^{-3}) \\
&= \frac{(k-1)(n-k) - k(n-k+1)}{nk} \binom{k}{2} w^2 \\
&\qquad + \frac{(k-1) w^2}{2k} + \frac{(k-1)^2 w^2}{2k} + O(k^2 n^{-3})
\end{aligned}
$$

and some algebraic manipulation reveals that all terms but the $O(k^2 n^{-3})$-term at the end cancel.

We have thus established for the bottom-to-top shuffle an eigenvalue of the form $\lambda = (1 - \gamma) e^{i\theta}$, where $\gamma = (1 + o(1)) \frac{2\pi^2 k(k-1)}{n^3}$ and $\theta = (1 + o(1)) \frac{2\pi}{n}$ and a corresponding eigenvector given by

$$
\mathbf{x} = [1, \lambda, \lambda^2, \ldots, \lambda^{n-k}, \lambda^{n-k}, \ldots, \lambda^{n-k}]^T.
$$

Now let $Z_t^j$ denote the position of card $j$ at time $t$, put $m = \lfloor n/2 \rfloor$, $\Phi^j(X_t) = x_{Z_t^j}$ and

$$
\Phi(X_t) = \sum_{j=1}^m \Phi^j(X_t).
$$

By linearity of expectation and what we just showed, $\mathbb{E}[\Phi(X_{t+1})|X_t] = \lambda\Phi(X_t)$. To apply Theorem 4.1 with $\Phi$ as test function, we need to check $\max_{s\in S}\mathbb{E}[|e^{-i\theta}\times\Phi(X_{t+1})-\Phi(X_t)|^2|X_t=s]$. However, deterministically, all the $n-k$ top cards move one step down the deck at each shuffle, and for such a card, $j$, on position $r$,

$$|e^{-i\theta}\Phi^j(X_{t+1})-\Phi^j(X_t)|=|e^{-i\theta}\lambda^{r+1}-\lambda^r|\leq\gamma=O(k^2n^{-3}).$$

For the card that is moved to the top, $\Phi^j$ changes from $\lambda^{n-k}$ to $1$, a change of order $O(kn^{-1})$. Finally, for $k-1$ of the cards at the bottom of the deck, their corresponding $\Phi^j$'s stay put at $\lambda^{n-k}$, so for these cards,

$$|e^{-i\theta}\Phi^j(X_{t+1})-\Phi_j(X_t)|=(1+o(1))|1-e^{iw}|=O(n^{-1}).$$

Using the triangle inequality,

$$|e^{-i\theta}\Phi(X_{t+1})-\Phi(X_t)|\leq mO(k^2n^{-3})+O(kn^{-1})+kO(n^{-1})=O(kn^{-1})$$

and, consequently,

$$\max_{s\in S}\mathbb{E}[|e^{-i\theta}\Phi(X_{t+1})-\Phi(X_t)|^2|X_t=s]=O(k^2n^{-2}).$$

We may thus apply Theorem 4.1 with $R=O(k^2n^{-2})$. Observing that $\Phi(s_0)=O(n)$ if we start with the cards in order (this is why we use only half the cards when defining $\Phi$), and putting $a=1/2$, we then get the following lower bound when $k=o(n)$:

$$\begin{aligned}
T&=(1+o(1))\frac{n^3}{2\pi^2k(k-1)}\left(\log\Phi(s_0)-\frac{1}{2}\log\frac{8R}{\gamma}\right)\\
&=(1+o(1))\frac{n^3}{2\pi^2k(k-1)}\left(\log n-\frac{1}{2}\log\frac{k^2n^{-2}}{k^2n^{-3}}\right)\\
&=(1+o(1))\frac{n^3}{4\pi^2k(k-1)}\log n.
\end{aligned}$$

Now what about the case $k=\Theta(n)$? In this case it is fairly easy to find a lower bound of order $n\log n$, just as one would expect from the $k=o(n)$ lower bound just given. This can be done either by using the $1-2/k$ eigenvalue and the corresponding eigenvector in Theorem 4.1 or by softer probabilistic reasoning as for the (inverse) random-to-top shuffle; after $\epsilon n\log n$ steps card, $n$ will still be very close to the bottom of the deck. However, an intriguing fact about these shuffles is that when $k$ is really large, for example, $k>2n/3$, then the second eigenvalue for the single card chain is farther away from $1$ than for the case $k=n$. This seems to suggest that the bottom-to-top shuffle with, say, $k=0.9n$, could in fact be faster then the random-to-top shuffle. Unfortunately we have not been able to find any good evidence about whether this is indeed the case or not, and it would be quite interesting to know the answer to this problem.

In summary, we have found the following result:

THEOREM 4.4. *A lower bound on the mixing time for the bottom-to-top shuffle where the card taken to the top at a given shuffle is chosen among the bottom $k = o(n)$ cards is given by*

$$(1 + o(1)) \frac{n^3}{4\pi^2 k(k-1)} \log n.$$

*When $k = \Theta(n)$, the mixing time is $\Omega(n \log n)$.*

Note that taking $k = 2$ in Theorem 4.4 recovers the previously known $\frac{1}{8\pi^2} n^3 \log n$ lower bound for the Rudvalis shuffle.

4.1.2. *GR shuffles with $p_{n-k} = p_n = 1/2$.* To avoid parity problems, we must assume that $k$ is odd. Using the same notation as for the bottom-to-top shuffle above, the equations for the eigenvalue/eigenvector pairs now become

$$x_1 = 1,$$
$$\lambda x_j = x_{j+1}, \qquad j = 1, 2, \ldots, n-k-1,$$
$$\lambda x_{n-k} = \tfrac{1}{2} + \tfrac{1}{2} x_{n-k+1},$$
$$\lambda x_j = \tfrac{1}{2} x_j + \tfrac{1}{2} x_{j+1}, \qquad j = n-k+1, n-k+2, \ldots, n-1,$$
$$\lambda x_n = \tfrac{1}{2} + \tfrac{1}{2} x_n.$$

Solving forward using all but the last equation, we get

$$x_1 = 1, \qquad x_2 = \lambda, \ldots, x_{n-k} = \lambda^{n-k-1},$$
$$x_{n-k+1} = 2\lambda^{n-k} - 1,$$
$$x_{n-k+2} = (2\lambda - 1)(2\lambda^{n-k} - 1), \ldots, x_n = (2\lambda - 1)^{k-1}(2\lambda^{n-k} - 1).$$

Invoking the last equation, we also have $x_n = 1/(2\lambda - 1)$ and since the two expressions for $x_n$ must agree, we get the characteristic equation

$$g(\lambda) := (2\lambda - 1)^k (2\lambda^{n-k} - 1) - 1 = 0.$$

Again we make separate treatments of the cases $k = o(n)$ and $k = \Theta(n)$, beginning with the former.

When $k = o(n)$, we suspect that a useful eigenvalue can be found quite close to 1. Therefore, $2\lambda - 1 \approx \lambda^2$, so it is a good idea to start with the equation one gets by making this replacement:

$$f(\lambda) := \lambda^{n+k} - \tfrac{1}{2}\lambda^{2k} - \tfrac{1}{2} = 0.$$

Putting $\mu := \lambda^k$ and $n = ck$, we get

$$\mu^{c+1} - \tfrac{1}{2}\mu^2 - \tfrac{1}{2} = 0.$$

Since $c = \Omega(1)$, this equation resembles the characteristic equation for the bottom-to-top shuffle. It has a solution very close to

$$\mu = \mu_0 := \left(1 - \frac{2\pi^2}{c^3}\right)e^{i2\pi/c}$$

and so a zero of $f(\lambda)$ can be found very close to

$$\lambda = \lambda_0 := \left(1 - \frac{k^2w^2}{2n}\right)e^{iw}.$$

Let us check how close to a solution $\lambda_0$ is:

$$\begin{aligned}
f(\lambda_0) &= \left(1 - \frac{k^2w^2}{2n}\right)^{n+k}e^{ikw} - \frac{1}{2}\left(1 - \frac{k^2w^2}{2n}\right)^{2k}e^{2ikw} - \frac{1}{2} \\
&= \left(1 - \frac{(n+k)k^2w^2}{2n} + O(k^4n^{-4})\right)\left(1 - \frac{k^2w^2}{2} + ikw + O(k^3n^{-3})\right) \\
&\quad - \frac{1}{2}(1 + O(k^3n^{-3}))(1 - 2k^2w^2 + 2ikw + O(k^3n^{-3})) - \frac{1}{2} \\
&= -\frac{n+k}{2n}k^2w^2 - \frac{1}{2}k^2w^2 + k^2w^2 + O(k^3n^{-3}) = O(k^3n^{-3}).
\end{aligned}$$

Since

$$f'(\lambda) = (n+k)\lambda^{n+k-1} - k\lambda^{2k-1} = (1 + o(1))n$$

in a large enough surrounding of $\lambda_0$, Lemma 4.2 implies that $f$ has a zero in a point $\lambda_1$ at distance $O(k^3n^{-4})$ from $\lambda_0$. Next we check how far off $\lambda_1$ is from being a zero of $g$. Compare $\lambda_1^2$ and $2\lambda_1 - 1$:

$$\begin{aligned}
2\lambda_1 - 1 - \lambda_1^2 &= 2\lambda_0 - 1 - \lambda_0^2 + O(k^3n^{-5}) \\
&= 2\left(1 - \frac{k^2w^2}{2n}\right)\left(1 - \frac{w^2}{2} + iw + O(n^{-3})\right) - 1 \\
&\quad - \left(1 - \frac{k^2w^2}{n} + O(k^4n^{-6})\right)(1 - 2w^2 + 2iw + O(n^{-3})) \\
&\quad + O(k^3n^{-5}) \\
&= w^2 + o(n^{-2}).
\end{aligned}$$

Thus,

$$(2\lambda_1 - 1)^k = \lambda_1^{2k} + kw^2 + o(kn^{-2}).$$

Therefore, since $\lambda_1$ is a zero of $f$,

$$\begin{aligned}
1 + g(\lambda_1) &= (\lambda_1^{2k} + kw^2 + o(kn^{-2}))(2\lambda_1^{n-k} - 1) \\
&= 1 + 2f(\lambda_1) + (1 + o(1))(kw^2 + o(kn^{-2})) \\
&= 1 + kw^2 + o(kn^{-2}),
\end{aligned}$$

so $g(\lambda_1) = (1 + o(1))kw^2$. However, $g'(\lambda) = (1 + o(1))2n$ in a sufficiently large surrounding of $\lambda_1$ to allow us to appeal to Lemma 4.2 for the conclusion that $g$ has a zero in a point $\lambda_2$ at distance $(1 + o(1))kw^2/(2n)$ from $\lambda_1$. By the triangle inequality,

$$|\lambda_2 - \lambda_0| = (1 + o(1))\frac{kw^2}{2n} + O(k^3 n^{-4})$$

and so $\lambda_2$, the eigenvalue we set out looking for, can be written as

$$\lambda = (1 - \gamma)e^{-i\theta},$$

where

$$(1 + o(1))\frac{2\pi^2 k(k-1)}{n^3} \le \gamma \le (1 + o(1))\frac{2\pi^2 k(k+1)}{n^3}$$

and $\theta = (1 + o(1))\frac{2\pi}{n}$.

Next we apply Theorem 4.1 with the eigenvalue and corresponding eigenvector, $\mathbf{x}$, we have just found. For card $j$, put $\Phi^j(X_t) = x_{Z_t^j}$, $m = \lfloor n/2 \rfloor$ and

$$\Phi(X_t) = \sum_{j=1}^{m} \Phi^j(X_t).$$

To bound $\Delta := |e^{-i\theta}\Phi(X_{t+1}) - \Phi(X_t)|$, observe that the top $n - k$ cards, just as for the bottom-to-top shuffle, deterministically move one step down the deck. Thus, each of them contribute $(1 + o(1))\gamma = O(k^2 n^{-3})$ to $\Delta$, so by the triangle inequality, they cannot together contribute more than $O(k^2 n^{-2})$. Of the $k$ cards between position $n - k$ and position $n$, $k - 1$ cards will move one step up the deck or stay put and will thereby individually contribute $O(n^{-1})$ and together at most $O(kn^{-1})$ to $\Delta$. Finally, one card will move from position one of the bottom $k$ positions to the top, thereby contributing $O(kn^{-1})$ to $\Delta$. Summing up, we get that $\Delta$ is bounded by $O(kn^{-1})$ and Theorem 4.1 may be applied with $R = O(k^2 n^{-2})$. Putting $a = 1/2$ and noting that $\Phi(s_0) = O(n)$ if $s_0$ has the cards in order, we get almost the same lower bound that we got for the bottom-to-top shuffles:

THEOREM 4.5. *Consider the GR shuffle where $p_{n-k} = p_n = 1/2$ and $k$ is odd. If $k = o(n)$, then a lower bound on the mixing time is given by*

$$(1 + o(1))\frac{n^3}{4\pi^2 k(k+1)}\log n.$$

Again note that we retrieve the previously known lower bound for the Rudvalis shuffle, this time by putting $k = 1$.

The case $k = \Theta(n)$ is again a little more difficult since it is harder to tell in general as exactly as for the case $k = o(n)$ where the second eigenvalue is to be found. However, it is fairly easy to show that there is an eigenvalue for which one has $\gamma = \Theta(n^{-1})$ and $\theta = O(n^{-1})$ and that there are no other nontrivial eigenvalues with $\gamma$ of lower order than this. Applying Theorem 4.1 with the corresponding eigenvector then gives a lower bound of $\Theta(n \log n)$, as expected from the lower bound for $k = o(n)$. With some care, one can come up with more explicit bounds if $k$ is specified, for example, when $k = n/2$, the same methods as those used above yield an eigenvalue $(1 - \gamma)e^{i\theta}$, where $\gamma = (1 + o(1))\frac{\log 2}{n}$ and $\theta = (1 + o(1))\frac{3}{4}w$. Applying Theorem 4.1 then gives the lower bound

$$T = (1 + o(1))\frac{1}{2\log 2}n \log n.$$

However, this is, at least not in general, the correct order of the mixing. To show this, let us give the case $k = n/2$ some extra consideration: Let $Z_t$ be the position of a specific card at time $t$ (counting $0, \ldots, n-1$ for the positions rather than $1, \ldots, n$). Put

$$U_t := Z_t + (Z_t - k)_+ - t \bmod k$$

and $V_t = U_t - U_{t-1} \bmod k$. Then $\mathbb{E}[V_t | X_{t-1}] = 0$ and $\mathrm{Var}(V_t | X_{t-1})$ is 1 when $U_{t-1} > k$ and 0 when $U_{t-1} \le k$. Thus,

$$\mathrm{Var}\, U_t = \mathbb{E}[\mathbb{E}[U_t^2 | X_{t-1}]] = \mathbb{E}[U_{t-1}^2 + \mathbb{E}[V_t^2 | X_{t-1}]] \le \mathrm{Var}\, U_{t-1} + 1,$$

so by induction, $\mathrm{Var}\, U_t \le t$. However, then, by Chebyshev's inequality,

$$P(|U_{10^{-6}n^2} - U_0 - 10^{-6}n^2| \ge 0.01n) \le 0.01,$$

while at stationarity a deviation like this would have a probability of more than 0.9. We have shown the following:

THEOREM 4.6.   *Suppose that $n$ is even. Then the mixing time of the GR shuffle with $p_n = p_{n/2} = 1/2$ is $\Omega(n^2)$.*

The proof of Theorem 4.6 can be made to work when $k = cn$ for any rational constant $c$ and then gives a lower bound of order $n^2$. However, we have not been able to turn this $\Omega(n^2)$ lower bound for single cards into an $\Omega(n^2 \log n)$ bound for the whole deck, mainly because of the strong dependence between the motions of different cards. When $c$ is irrational, things are more unclear; it is not hard to see that in this case the mixing time for single cards is in fact $\Theta(n \log n)$, but it seems hard to imagine that the whole deck would mix much faster for irrational values of $c$.

What the true order of mixing is for $p_n = p_{n-k} = 1/2$ seems to be a wide open question. One natural guess is that the mixing time is $\Theta(n^3 \log n)$ for all $k$ and either a proof or a counterexample to this would be very interesting.

4.2. *Upper bounds.* As already pointed out, we will only be able to provide a good upper bound for the bottom-to-top shuffle. Thus, we are again considering the situation where at each time step the card moved to the top of the deck is chosen uniformly at random from the $k$ bottom cards. As usual, denote the state of the deck at time $t$ by $X_t$, and let $X_0$ be any fixed state. (Since the bottom-to-top shuffle describes a simple random walk on $S_n$, the stationary distribution is uniform and the particular starting state does not affect the convergence rate.) For each $t = 0, 1, 2, 3, \ldots$, let the mapping $\beta_t : S_n \to S_n$ be given by

$$\beta_t(\sigma) = (n-1\ n-2 \ldots 1\ 0)^t \circ \sigma$$

(where the positions of the permutations are for convenience now denoted $0, 1, 2, \ldots, n-1$). Put $Y_t = \beta_t(X_t)$. In words, $Y_t$ is $X_t$ with position $k+t$ (modulo $n$) regarded as position $k$, $k = 0, 1, \ldots, n-1$. Clearly,

$$\|P(Y_t \in \cdot) - \pi\| = \|P(X_t \in \cdot) - \pi\|,$$

so we may, and will, work with $Y_t$ instead of $X_t$ itself. The process $\{Y_t\}$ describes the (time-inhomogenous) Markov chain one gets by at time $t$, making a random-to-bottom shuffle modulo $n$ on the subset of positions $A_t := \{n-k+1-t, n-k+2-t, \ldots, n-t\}$ modulo $n$. (The set $A_t$ corresponds to the bottom $k$ positions for $\{X_t\}$.)

We will couple another process $\{Y_t'\}$ with the same updating mechanism, but started from stationarity, with $\{Y_t\}$. The coupling rule is the following: For all cards, $c$, such that $Y_t(c) \in A_t$ and $Y_t'(c) \in A_t$, move $c$ to position $n-t$ in $Y_t'$ if and only if it is moved also in $Y_t$. If the card moved in $Y_t$ is not in an $A_t$-position in $Y_t'$, then pick for $Y_t'$ a card chosen uniformly from those cards that are in an $A_t$-position in $Y_t'$ but not in $Y_t$. Clearly, $\{Y_t'\}$ has the correct updating distribution and the coupling has the following important properties:

- A card $c$ in $\{Y_t\}$ cannot pass its copy in $\{Y_t'\}$ unless $Y_t(c)$ and $Y_t'(c)$ are both in $A_t$.
- Once a card $c$ has been in an $A_t$-position in $Y_t$ and $Y_t'$ simultaneously, it will be matched as soon as it leaves $A_t$. This is because the only way a card can leave $A_t$ is by being picked as the card moved to the bottom of $A_t$ and by the nature of the coupling, this happens simultaneously for the two decks.

Consequently, let $T_0(c)$ be the first time that card $c$ has been in $A_t$ simultaneously for the two decks:

$$T_0(c) := \min\{t : Y_t(c) \in A_t \text{ and } Y_t'(c) \in A_t\}$$

and let $T_0 = \max_c T_0(c)$. Then at time $T_0$ all cards outside $A_{T_0}$ must be matched. However, then, by the nature of the shuffle and the coupling, *all* cards will be matched as soon as all cards in $A_{T_0}$ have left $A_t$. Putting $T$ for the first time this has happened, the coupling inequality tells us that

$$\|P(Y_t \in \cdot) - \pi\| \le P(T > t),$$

so for the rest of this section, we focus on estimating $T$. By the standard coupon collector's problem, for any $a > 0$, $P(T - T_0 > (1+a)k \log k) = o(1)$ unless $k = O(1)$, in which case $P(T - T_0 > f_n) = o(1)$ as soon as $f_n = \Omega(1)$. Now what about $T_0$?

Consider the motion of a single card $c$ in $\{Y_t\}$. Denote the time interval between two successive occasions when $c$ leaves $A_t$, by a *cycle* for $c$. A moment's thought reveals that during a cycle $c$ moves $k - G$ steps down the deck (modulo $n$), where $G$ is a random variable with geometric distribution with parameter $1/k$. Thus, putting $\tau_j = \tau_j(c)$ for the $j$th time, $j = 1, 2, 3, \ldots$, that $c$ leaves $A_t$, the process $\{Y_{\tau_j}(c)\}_{j=1}^{\infty}$ describes a random walk on $\mathbb{Z}_n$ with step size distribution $k - G$. Keeping track of how $c$ crosses the top/bottom border of the deck, we get a random walk on $\mathbb{Z}$ with step size distribution $k - G$; in particular, the step size mean is 0 and the step size variance is $k^2(1 - 1/k) = k(k-1)$. Now considering $c$'s motion in $\{Y_t'\}$ gives a corresponding sequence $\{\tau_j'\}$ of stopping times and a corresponding random walk with the same properties. Let $J = \min\{j : \tau_j' > T_0(c)\}$. Then if $Y_0(c) > Y_0'(c)$,

$$\tau_1 < \tau_1' < \tau_2 < \tau_2' < \cdots < \tau_J = \tau_J' < \tau_{J+1} = \tau_{J+1}' < \cdots$$

or

$$\tau_1 < \tau_1' < \tau_2 < \tau_2' < \cdots < \tau_{J-1}' = \tau_J < \tau_J' = \tau_{J+1} < \cdots,$$

depending on if $T_0(c)$ coincides with $c$ entering $A_t$ in $Y_t'$ or $Y_t$. If $Y_0(c) < Y_0'(c)$, then these relations with the primed and nonprimed quantities interchanged hold. In the sequel we assume that $Y_0(c) > Y_0'(c)$; the other case is treated analogously.

Put $S_j = Y_{\tau_j}(c) - Y_{\tau_j'}'(c)$. Then $\{S_j\}$ for the first $J - 1$ steps behaves like the difference between two independent random walks on $\mathbb{Z}_n$. Thus, $\{S_j\}$ itself for the first $J - 1$ steps describes a symmetric random walk on $\mathbb{Z}$ with step size variance $2k(k-1)$. To make this exact, define a third deck $\{Y_t''\}$ such that $Y_t''$ coincides with $Y_t'$ for $t < T_0(c)$, but let $Y_t''$ evolve independently of $Y_t$ from time $T_0(c)$ and on. Associate in analogy with the above with $Y_t''$ the stopping times $\tau_j''$ and put $U_j = Y_{\tau_j}(c) - Y_{\tau_j''}''(c)$. Then $\{U_j\}$ describes a random walk, started from somewhere between 0 and $n$, that for all time is the difference between two random walks of the kind encountered above, and $U_j$ coincides with $S_j$ for the first $J - 1$ steps. The process $\{U_j\}$ also has the property that if it, after $j$ steps, has at least once passed the origin or vertex $n$, then $j \ge J$; this follows from the above properties of the coupling. Thus, $T_0(c)$ is bounded by the first time $U_j$ passes outside the interval $[0, n]$.

Let us now bound the probability that $\{U_j\}$ has not passed outside $[0, n]$ in, say, $j_0$ steps. This probability is maximized when $U_0 = n/2$, so let us assume that this is the case. Put

$$W_j := \frac{1}{\sqrt{2k(k-1)}}\left(U_j - \frac{n}{2}\right),$$

so that $W_0 = 0$, $\{W_j\}$ has step size variance 1 and $U_j$ passes outside $[1, n-1]$ when $W_j$ passes outside

$$\left[-\frac{n}{2^{3/2}\sqrt{k(k-1)}}, \frac{n}{2^{3/2}\sqrt{k(k-1)}}\right].$$

From here we must treat the cases $k = o(n)$ and $k = \Theta(n)$ separately. Assume first that $k = o(n)$. Let

$$M := \frac{n}{2^{3/2}\sqrt{k(k-1)}}.$$

By Donsker's theorem (since $M \to \infty$) (see, e.g., [3], Section 7.6)

$$\left\{\frac{1}{M}W_{M^2 s}\right\}_{s \in [0, \infty)} \xrightarrow{D} \{B_s\}_{s \in [0, \infty)}$$

as $n \to \infty$, where $\{B_s\}$ stands for standard Brownian motion. Thus,

$$P\left(\forall s \le s_0 : \left|\frac{1}{M}W_{M^2 s}\right| < 1\right) = (1 + o(1))P(\forall s \le s_0 : |B_s| < 1)$$

$$\le (1 + o(1))\frac{4}{\pi}e^{-\pi^2 s_0/8},$$

where the last bound can be found, for example, in [3], Section 7.8. With $s_0 = (1 + o(1))(8/\pi^2)\log n$, the right-hand side is $o(n^{-1})$. Translating this back to $\{U_j\}$ tells us that the probability that $U_j$ has not passed outside $[0, n]$ after $(1 + o(1))(8/\pi^2)M^2 \log n$ steps is $o(n^{-1})$. Thus, with

$$j_0 := (1 + o(1))\frac{n^2}{\pi^2 k(k-1)}\log n,$$

we have

$$P(J > j_0) = o(n^{-1}).$$

Now how long does it take before $c$ has gone through $j_0$ cycles? This time, say, $T_1$, can be written as

$$T_1 = \eta + \sum_{j=1}^{j_0} \xi_j,$$

where $\eta$ is the time taken until $c$ leaves $A_t$ for the first time, and the $\xi_j$'s are the $j_0$ independent cycle times. Since $\eta$ represents the time taken for a partial cycle, $\eta$ is stochastically dominated by a random variable $\xi_0$ with the same distribution as $\xi_1, \ldots, \xi_{j_0}$ and so $T_1$ is dominated by $T_2 := \sum_{j=0}^{j_0} \xi_j$. The $\xi_j$'s have distribution $n - k + G$, where $G$ is geometric with parameter $1/k$. Hence, $\mathbb{E}T_2 = n(j_0 + 1)$ and we want to bound $P(T_2 \geq (1+a)n(j_0 + 1))$ for an arbitrary small $a > 0$. However, this probability coincides with the probability that $B < j_0$, where $B$ is binomial with parameters $(an + k)(j_0 + 1)$, and $1/k$. Since $j_0 = \Omega(\log n)$ and $k = o(n)$, it follows from standard Chernoff bounds that

$$P(T_2 > (1+a)nj_0) = o(n^{-1}).$$

Thus, $P(\tau'_{j_0} > (1+a)nj_0) = o(n^{-1})$ and since it was shown above that $P(J > j_0) = o(n^{-1})$ and since $T_0(c) < \tau'_J$, we get that

$$P(T_0(c) > (1+a)nj_0) = o(n^{-1}).$$

Summing over the cards,

$$P(T_0 > (1+a)nj_0) = o(1).$$

Finally, since $k = o(n)$, $k \log k$ is of smaller order than $nj_0$ [or when $k = O(1)$, then $f_n$ is of smaller order than $nj_0$], and so

$$P(T > (1+a)nj_0) = o(1).$$

We have thus arrived at the following upper bound on the mixing time:

$$\tau_{\mathrm{mix}} \leq (1 + o(1)) \frac{n^3}{\pi^2 k(k-1)} \log n.$$

In the case $k = \Theta(n)$, the above approach does not work exactly as it stands for two reasons:

- Since $n$ and $k$ are of the same order, the Brownian motion approximation does not work.
- The $k \log k$ term for $T - T_0$ is of the same order as $nj_0$.

However, it is easy to modify the arguments slightly to give an upper bound of the type $Cn \log n$ for some constant $C$. A more detailed analysis would also give an estimate for $C$, but since these estimates appear to be well above what one would guess are the correct figures (e.g., when $k$ is close to $n$, an estimate of $C$ lands well above 1), we omit that here.

The following theorem summarizes our results on the bottom-to-top shuffle:

THEOREM 4.7.  *Let $\tau_{\mathrm{mix}}$ be the mixing time for the bottom-to-top shuffle where the card taken to the top at each shuffle is chosen among the bottom $k$ cards. Then if $k = o(n)$,*

$$(1 + o(1)) \frac{n^3}{4\pi^2 k(k-1)} \log n < \tau_{\mathrm{mix}} < (1 + o(1)) \frac{n^3}{\pi^2 k(k-1)} \log n.$$

*If $k = \Theta(n)$, then $\tau_{\mathrm{mix}} = \Theta(n \log n)$.*

**Acknowledgments.** I would like to thank Laurent Saloff-Coste and Yuval Peres for valuable discussions. I am also grateful to Sven Järner and Johan Karlsson for their help on Lemma 4.2. Finally, I would like to thank the referee for a thorough reading and many valuable comments.

## REFERENCES

[1] ALDOUS, D. and DIACONIS, P. (1986). Shuffling cards and stopping times. *Amer. Math. Monthly* **93** 333–348. MR0841111

[2] BAYER, D. and DIACONIS, P. (1992). Trailing the dovetail shuffle to its lair. *Ann. Appl. Probab.* **2** 294–313. MR1161056

[3] DURRETT, R. (1991). *Probability: Theory and Examples.* Wadsworth and Brooks/Cole, Pacific Grove, CA. MR1068527

[4] GOEL, S. (2006). Analysis of top to bottom-$k$ shuffles. *Ann. Appl. Probab.* **16** 30–55.

[5] HILDEBRAND, M. (1990). Rates of convergence of some random processes on finite groups. Ph.D. dissertation, Harvard Univ.

[6] JONASSON, J. (2006). The overhand shuffle mixes in $\Theta(n^2 \log n)$ steps. *Ann. Appl. Probab.* **16** 231–243.

[7] LINDVALL, T. (1992). *Lectures on the Coupling Method.* Wiley, New York. MR1180522

[8] MOSSEL, E., PERES, Y. and SINCLAIR, A. (2004). Shuffling by semi-random transpositions. In *Proceedings of the 45th Annual IEEE Symposium on Foundations of Computer Science* 572–581. IEEE, Los Alamitos, CA.

[9] RODRIGUES, E. (1995). Convergence to stationarity for a Markov move-to-front scheme. *J. Appl. Probab.* **32** 768–776. MR1344075

[10] SALOFF-COSTE, L. (2004). Total variation lower bounds for Markov chains: Wilson's lemma. In *Random Walks and Geometry* (V. A. Kaimanovich, K. Schmidt and W. Goess, eds.) 515–532. de Gruyter, Berlin. MR2087800

[11] WILSON, D. B. (2003). Mixing time of the Rudvalis shuffle. *Electron. Comm. Probab.* **8** 77–85. MR1987096

[12] WILSON, D. B. (2004). Mixing times of lozenge tiling and card shuffling Markov chains. *Ann. Appl. Probab.* **14** 274–325. MR2023023

DEPARTMENT OF MATHEMATICS
CHALMERS UNIVERSITY OF TECHNOLOGY
S-412 96 GÖTEBORG
SWEDEN
E-MAIL: jonasson@math.chalmers.se